# Integer conversions and estimation of the number of integer solutions for algebraic Diophantine equations with nondiagonal form

VICTOR VOLFSON

ABSTRACT. The paper assesses the top number of integer solutions for algebraic Diophantine Thue diagonal equation of the degree $n \geq 2$ and number of variables $k > 2$ and equations with explicit variable in the case when the coefficients of the equation are of the opposite signs. The author found integer conversions that maintain the asymptotic behavior of the number of integer solutions of algebraic Diophantine equation in the case of the conversion equation to diagonal form. The paper considers the estimation of the number of integer solutions for some types of algebraic Diophantine equations with nondiagonal form.

1. INTRODUCTION

Algebraic Diophantine equation of the degree n, variables k and integer coefficients has the form: $\sum_{i=0}^{n} F_i(x_1,...x_k) = 0$, where $F_i(x_1,...x_k)$ is a form for variables k with degree i and $F_0$ is an integer.

Thue showed that an irreducible equation $F_n(x_1,...x_k) + F_0 = 0$ ($n \geq 3$, $k = 2$ and $F_0 \neq 0$) has only a finite number of integer solutions and in the particular case it can have no solution at all. Thue method was developed in Siegel's works [1], which set well-known theorem about the finiteness of the number of integer points on a curve of genus $g \geq 1$.

This method was extended by Schmidt [2] to the case $k > 2$, which allowed him to receive a multi-dimensional generalization of Thue's results of a finite number of integer solutions for the normative Diophantine equation.

___

Currently, quantitative estimates for the integer solutions of some classical Diophantine equations obtained by this method. Baker [3] made an effective assessment of the Tue equation. A similar effective analysis for more general Tue-Mahler equation: $f(x,y) = mp_1^{x_1}...p_s^{x_s}$ (where $p_1,...,p_s$ are fixed numbers) was carried out in [4]. Another class of Diophantine equations (that allow effective analysis) constitutes superelliptic equation: $y^s = f(x)$, where $s \geq 2$, $f$ is an integer function of degree $n \geq 3$. Effective analysis of this equation was conducted by Baker [5]. This result was significantly strengthened [6]. Thus, the effective analysis is done only for algebraic Diophantine equations of two variables

The Hardy and Littlewood circular method (CM) [7] allows us to make an upper estimate of the number of natural solutions to various algebraic Diophantine diagonal equations with integer coefficients and large number of variables. The method has been significantly strengthened Vinogradov [8]. An upper estimate of the number of natural solutions was obtained in [9] (with the help of the CM) of the diagonal Thue equation: $a_1 x_1^n + ... + a_k^n + a_0 = 0$ (in the case when all coefficients $a_i$ are integers and of the same sign, and $a_0$ is integer (not equal to 0) and of another sign.

CM was used in [10] to obtain upper estimations of the number of natural solutions of various homogeneous algebraic diagonal Diophantine equations with integer coefficients ($a_0 = 0$).

The number of natural solutions of inhomogeneous algebraic diagonal Diophantine equation, including diagonal Thue equation with coefficients of the same sign. was estimated in [11]

It is interesting to find estimates of the number of integer (natural) solutions of Thue's equation for values $k > 2$, $n \geq 2$ when coefficients of Thue's equation have different signs. It is also interesting the consideration of the inhomogeneous diagonal Diophantine equation with explicit variable. These issues are considered in the section 2 of the paper.

Estimates of the number of integer solutions of Diophantine equations in [9], [10], [11] were carried out mainly for the diagonal equations, which correspond to the canonical equations of surfaces. Therefore, there is an important question could any algebraic Diophantine equation be reduced to diagonal form? The answer to this question is in the section 3 of the paper.

## 2. THE ESTIMATION OF THE NUMBER OF INTEGER (NATURAL) SOLUTIONS FOR SOME TYPES OF ALGEBRAIC DIAGONAL DIOPHANTINE EQUATIONS

Let us consider the inhomogeneous Diophantine diagonal equation with explicit variable:

$$x_1 = b_2 x_2^n + \ldots + b_k x_k^n, \qquad (2.1)$$

where $b_i$ are integers, and $k$, $n$ are natural numbers.

The case when all $b_i$ are natural is shown in [11]. Now we'll examine the case when values $b_i$ have different signs. We get the equation of a hyperbolic paraboloid (values $k = 3$ and $n = 2$):

$$x_1 = a_2 x_2^2 - a_3 x_3^2, \qquad (2.2)$$

where $a_2, a_3$ are natural numbers.

Let us consider cross-section of a hyperbolic paraboloid (2.2) with the plane $x_3 = h$ and we obtain parabola shifted from the coordinate origin:

$$x_1 + a_3 h^2 = a_2 x_2^2. \qquad (2.3)$$

It is well known that we have the following asymptotic estimate for the number of integer points for the parabola (2.3) in the square with the side $[-N, N]$: $R_2(N) \ll N^{1/2}$.

The hyperbolic paraboloid (2.2) is not limited by the axis $x_3$, so the segment of the axis $[-N, N]$ has $2N + 1$ such parabolas or $O(N)$.

Therefore the asymptotic upper bound for the number of integer solutions of the equation (2.2) (in the cube with the sides $[-N, N]$) is equal to:

$$R_3(N) \ll N^{3/2}. \qquad (2.4)$$

The equation (2.1) can have the following form:

$$x_1 = a_2 x_2^n + \ldots + a_m x_m^n - a_{m+1} x_{m+1}^n - \ldots - a_k x_k^n, \qquad (2.5)$$

where all values $a_i$ are natural numbers.

We show (by induction method) that the upper asymptotic estimate of the number of integer solutions of the equation (2.5) in the hupercube with the sides $[-N, N]$) is equal to:

$$R_k(N) \ll N^{(m-1)/n+k-m}. \tag{2.6}$$

We consider the equation in the first step of induction:

$$x_1 = a_2 x_2^n + ... + a_m x_m^n - a_{m+1} x_{m+1}^n. \tag{2.7}$$

If crossing the hypersurface (2.7) by the plane $x_{m+1} = h$, we get the hypersurface shifted from the coordinate origin:

$$x_1 + a_{m+1} h^n = a_2 x_2^n + ... + a_m x_m^n. \tag{2.8}$$

As shown in [11] we have the following asymptotic estimation for the number of integer solutions of the equation (2.8) in a hypercube with the sides $[-N, N]$:

$$R_m(N) \ll N^{(m-1)/n}. \tag{2.9}$$

Since the hypersurface (2.7) is not limited in axis $x_{m+1}$, it comprises $2N+1$ or $O(N)$ sections (2.8) on the axis of the segment $[-N, N]$. Therefore, we have the following asymptotic estimate for the number of integer solutions for the equation (2.7) in a hypercube with the sides $[-N, N]$:

$$R_{m+1}(N) \ll N^{(m-1)/n+1}. \tag{2.10}$$

Now suppose that (2.6) is right for $k = l$:

$$R_l(N) \ll N^{(m-1)/n+l-m}. \tag{2.11}$$

is taken place.

Having in mind (2.11), we'll prove that the estimate (2.6) is right for the value $k = l+1$:

$$R_{l+1}(N) \ll N^{(m-1)/n+l-m} \cdot N = N^{(m-1)/n+l+1-m},$$

since there is $O(N)$ cross section with the estimate (2.11) QED

Formula (2.4) is obtained from (2.6) for values: $n = 2, k = 3, m = 2$.

Formula (2.6) is right for even values of $n$, as well as for non-negative values of the variables for an odd value of $n$.

We use Pila formula [12] for estimating an upper estimate for the number of integer solutions of irreducible equations of the degree $n$ in the hypercube with a side $[-N, N]$, which is easily proved by induction of the number of variables $k$:

$$R_k(N) << N^{k-1+1/n+\epsilon}. \tag{.2.12}$$

Now we consider an upper estimate for the number of integer solutions of Thue equation. It is well known that the number of integer solutions for Thue equation with values $k=2, n>2$ is finite. Here we consider other cases.

First, we consider the diagonal Thue equation wih the number of variables $k>2$ and of the degree $n=2$:

$$b_1 x_1^2 + ... + b_k x_k^2 + c = 0. \tag{2.13}$$

The equation (2.13) corresponds to the imaginary ellipsoid and has no real solutions, if $b_i, c$ are the same sign. The equation (2.13) corresponds to an ellipsoid and has a finite number of integer solutions, if all $b_i$ are the same sign, and $c$ is the opposite sign. Therefore, the following estimate is true:

$$R_k(N) << N^{\epsilon}, \tag{2.14}$$

where $\epsilon$ is a small positive number.

The equation (2.13) can have the following form (if $b_i$ are of the sign opposite to the sign $c$ when $i \leq m$ and if $b_i$ and $c$ are the same sign when $i > m$):

$$a_1 x_1^2 + ... + a_m x_m^2 - a_{m+1} x_{m+1}^2 - ... - a_k x_k^2 + c = 0, \tag{2.15}$$

where values $a_i$ are natural numbers, and $c$ is an integer number.

The equation (2.15) corresponds to the hyperboloid with the index $m$. We have the following upper bound for the number of integer solutions of the equation (2.15) in the hypercube with the sides $[-N, N]$:

$$R_k(N) << N^{k-m+\epsilon}. \tag{2.16}$$

We'll prove (2.16) in a more general case a little later.

Now we note that when the size $N$ tends to infinity, the hyperboloid (2.15) tends to its asymptotic cone of directions:

$$a_1 x_1^2 + ... + a_m x_m^2 - a_{m+1} x_{m+1}^2 - ... - a_k x_k^2 = 0, \tag{2.17}$$

i.e. the homogeneous equation (2.17) has the same asymptotic behavior of the number of integer solutions as the equation (2.16).

Now we consider a more general diagonal Thue equation:

$$a_1 x_1^n + ... a_m x_m^n - a_{m+1} x_{m+1}^n - ... - a_k x_k^n + c = 0, \tag{2.18}$$

where all values $a_i$, $n$ are natural numbers, and $c$ is an integer one.

Using the method of mathematical induction we'll show that the number of integer solutions for the equation (2.18) has the estimate (2.16) if $n$ is of even number in the hypercube with the sides $[-N, N]$.

In the first step, we consider the equation:

$$a_1 x_1^n + ... a_m x_m^n - a_{m+1} x_{m+1}^n + c = 0. \tag{2.19}$$

If crossing the hypersurface (2.19) wit the hyperplane $x_{m+1} = h$, then we obtain the equation: $a_1 x_1^n + ... a_m x_m^n = a_{m+1} h^n - c$, which has a finite number of solutions i.e. $R_m(N) \ll N^\epsilon$, where $\epsilon$ is a small positive number. This corresponds to equation (2.15) for the value of $k = m$.

Let us assume:

$$R_l(N) \ll N^{l-m+\epsilon}. \tag{2.20}$$

is performed for the value $k = l$.

Let the hyperplane $x_{l+1} = h$ intersect the hypersurface:

$$a_1 x_1^n + ... a_m x_m^n - a_{m+1} x_{m+1}^n - ... - a_{l+1} x_{l+1}^n + c = 0. \tag{2.21}$$

We get the following hypersurface in the section:

$$a_1 x_1^n + \ldots a_m x_m^n - a_{m+1} x_{m+1}^n - \ldots - a_l x_l^n = a_{l+1} h_{l+1}^n - c.$$

The estimate for its number of integer solutions is (2.20).

Since there are no more than $2N+1$ or $O(N)$ of such sections, and based on (2.20) we obtain an estimate for the number of integer solutions for the equation (2.21): $R_{l+1} \ll N^{l-m+\epsilon} \cdot N = N^{l-m+1+\epsilon}$, that corresponds to (2.16).

Formula (2.16) is right for estimating the number of natural solutions for the equation (2.18) for the hypercube with the sides $[-N, N]$ if value of $n$ is odd.

It is necessary to use Pila formula (2.12) for estimating the number of top all integer solutions for the equation (2.18) if an odd value of n is in the hypercube with the sides $[-N, N]$.

We have the upper bound (2.16) for the number of integer solutions of equation (2.22) in the hypercube with the side $[-N, N]$ as the asymptotic cone of directions for the hypersurface (2.18) is the hypersurface:

$$a_1 x_1^n + \ldots a_m x_m^n - a_{m+1} x_{m+1}^n - \ldots - a_k x_k^n = 0, \qquad (2.22)$$

If the equation (2.18) has more negative factors than positive ones ($k - m > m$ or $k > 2m$), then the equation can be replaced by the equivalent one:

$$a_{m+1} x_{m+1}^n + \ldots + a_k x_k^n - a_1 x_1^n - \ldots - a_m x_m^n - c = 0,$$

with more positive factors. Therefore the upper bound of the number of integer solutions for the equation (2.18) in the hypercube with sides $[-N, N]$ does not exceed

$$R_k(N) \ll N^{[k/2]+\epsilon}. \qquad (2.23)$$

Formula (2.16) is right for small and large values of $k$. Note that the previously we considered assessment for the formula with the number of integer solutions of the homogeneous diagonal equation (of the degree $n$) for large number of variables:

$$a_1 x_1^n + \ldots + a_k x_k^n = 0, \qquad (2.24)$$

where all values $a_i$ are integers.

Based on Hardy and Littlewood circle method the formula for the top estimating of the number of integer solutions for the equation (2.24) with a large number of variables $k$ in the hypercube with the sides $[-N, N]$ was proved:

$$R_k \ll N^{k-n+\epsilon}, \tag{2.25}$$

where $\epsilon$ is a small positive number.

Formula (2.25) is right for Thue equation with a large number of variables:

$$a_1 x_1^n + ... + a_k x_k^n + c = 0, \tag{2.26}$$

since (2.24) is the equation of the asymptotic cone of directions.

Now we consider the estimation of the number of integer solutions for the equation:

$$x_1^2 + x_2^2 - x_3^2 - 1 = 0. \tag{2.27}$$

The equation (2.27) is an equivalent to the equation:

$$(x_1 + x_3)(x_1 - x_3) + (x_2 + 1)(x_2 - 1) = 0. \tag{2.27}$$

with integer solutions: $x_1 = -x_3, x_2 = 1$; $x_1 = x_3, x_2 = 1$; $x_1 = x_3, x_2 = -1$; $x_1 = -x_3, x_2 = -1$.

Therefore, this equation has the following number of integer solutions in the cube with the sides $[-N, N]$:

$$R_3(N) = 4(2N + 1) = O(N). \tag{2.28}$$

On the other hand, if having in mind (2.16), the estimate for the equation (2.27) is $R_3(N) \ll N^{1+\epsilon}$, that corresponds to (2.28).

3. INTEGER TRANSFORMATIONS OF DIOFHANTINE EQUATIONS WITH NONDIAGONAL FORM

Up to now we considered the asymptotic estimate of the number of integer solutions for algebraic Diophantine diagonal equations. Therefore, we are interested in integer transform, which would result the equation (3.1) to diagonal form, and thus would maintain the asymptotic behavior of the number of integer solutions for the equation:

$$F(x_1,...x_n) = 0. \tag{3.1}$$

Thus, we are interested in integer transform, which keeps the number of integer solutions for the equation, multiplied by real constant.

In this case each integer solution of the equation (3.1) must correspond to an integer solution of the equation $G(x^1_1,...x^1_n) = 0$ and vice versa. The polynomial $G(x^1_1,...x^1_n)$ is obtained after transformation of the polynomial $F(x_1,...x_n)$. Thus, the transformation must be an integer and bijective.

However, not every bijection mapping maintains the number of solutions, for example, for a hypersphere.

It is known that an algebraic equation corresponds to an algebraic surfaces and a diagonal algebraic equation corresponds to the canonical form of an algebraic surface. Algebraic surfaces of the second order can get a canonical form by a specified conversion movement. In some particular cases it can be done for surfaces of higher order [13]. In general, it cannot be done for any algebraic surface.

A conversion movement is a bijective mapping; it maintains the number of solutions for the hypersphere because it maintains the distance between points, for example from the point of solving the equation to the coordinate origin.

An integer conversion movement, which results in the canonical form of an algebraic surface, consists of the rotation converting matrix which members are only of -1, 0, 1, and integer transfer.

However, we don'tt need to keep the number of integer solutions of the equation (3.1) for a hypersphere. As mentioned above it is sufficient if using integer transform to keep the number of integer solutions, multiplied by a real constant.

Thus, one of the integer transformations maintaining the asymptotic behavior for the number of solutions of algebraic Diophantine equation is the deformation, which is of no integer motion transforming, i.e. surface reduction to the canonical form and performing a homothety with real constant k.

A purpose is to show when it is possible and how to choose the value of k in these cases.

In general the conversion movement that makes it possible to get the canonical form of the surface has the form:

$$x_1 = c_{11}x'_1 + \ldots + c_{1n}x'_n + c_1, \ldots, \quad x_n = c_{n1}x'_1 + \ldots + c_{nn}x'_n + c_n, \qquad (3.2)$$

where $c_{ij}$ $(i=1,\ldots,n; j=1,\ldots,n)$ are the rotation matrix coefficients, and $c_i$ $(i=1,\ldots,n)$ are the coordinates of the parallel transport.

Total transformation after homothety of $x'_i = kx''_i$ has the form:

$$x_1 = c_{11}kx''_1 + \ldots + c_{1n}kx''_n + c_1, \ldots, \quad x_n = c_{n1}kx''_1 + \ldots + c_{nn}kx''_n + c_n. \qquad (3.3)$$

The transformation (3.3) is an integer, if two conditions are met:

1. Transfer Coordinates of $c_i$ are integers.

2. All transform coefficients $kc_{ij}$ are integers.

The first point answers the question when such an integer transform is possible. Assume that it is taken place.

The second point answers the question what the value of $k$ should be to ensure all transform coefficients to be integers.

Let us consider the following cases of values of rotation matrix coefficients.

1. Any rotation matrix coefficients are correct rational fractions $c_{ij} = p_{ij}/q_{ij}$ or integers (-1,0,1), then it requires

$$k = LCM_{ij}(q_{ij}), \qquad (3.4)$$

for the implementation of the point 2 where $k = LCM$ is the least common multiple of all $q_{ij}$.

2. All rotation matrix coefficients contain the same irrationality $\sqrt{m}$. In this case it is required:

$$k = \sqrt{m}. \qquad (3.5)$$

to implement the point 2.

3. Various rotation matrix coefficients contain $k$ different irrationals. In this case, such $k$ does not exist.

4. Some rotation matrix coefficients are rational fractions and some factors contain irrationality. In this case $k$ also does not exist.

Let us consider the first case. We'll show that there is the case where all of the rotation matrix coefficients are rational fractions $p_i/q_i \leq 1$ and find its analytical expression for the value $n = 2$.

The proper rotation matrix looks

$$C = \begin{pmatrix} p_1/q_1 & -p_2/q_2 \\ p_2/q_2 & p_1/q_1 \end{pmatrix}, \qquad (3.6)$$

for value $n = 2$, where

$$(p_1/q_1)^2 + (p_2/q_2)^2 = 1. \qquad (3.7)$$

Equation (3.7) is Fermat's equation for the value $n = 2$:

$$(p_1 q_2)^2 + (p_2 q_1)^2 = (q_1 q_2)^2. \qquad (3.8)$$

Therefore, we have the following expression for the matrix coefficients in (3.6):

$$p_1/q_1 = p_1 q_2 / q_1 q_2 = (m^2 - n^2)/(m^2 + n^2), \quad p_2/q_2 = p_2 q_1 / q_1 q_2 = (2mn)/(m^2 + n^2), \quad (3.9)$$

where $m, n$ are natural numbers and $m > n$.

Transformations own rotation of C for $n = k (k > 2)$) consists of $k$ consecutive rotations around the coordinate axes. Rotation matrix $C_i$ around the i-th axis contains the value 1 at the intersection of the $i$-th row and $i$-th column and there are rotation matrix coefficients of $k-1$ order at the remaining rows and columns. Therefore, there is formula:

$$C = \prod_{i=1}^{k} C_i. \qquad (3.9)$$

Now we'll prove by induction that there is a proper rotation matrix of any order, which members are the only rational fractions $p_{ij}/q_{ij} \leq 1$.

We have already shown this for the proper rotation matrix of the order n = 2. Let us suppose that there is the proper rotation matrix C of order $n = k$ which members are only rational fractions $p_{ij}/q_{ij} \leq 1$. Then the proper rotation matrix C of the order $n = k+1$ (because of (3.9)) is the product of matrices $C_i$. The value 1 and regular rational fractions are members of the matrices $C_i$ of order $n = k$. Therefore, members of the proper rotation matrix of order $n = k+1$ are the sum of products of 1 and regular rational fractions ie they are rational fractions. Since the matrix C is a proper rotation matrix, then the following condition $p_{ij}/q_{ij} \leq 1$ is satisfied for its members QED.

Now we show that there is a rotation matrix C for the second case.

Let us consider the homogeneous Diophantine equation with two variables:

$$a_{11}x_1^2 + 2a_{12}x_1x_2 + a_{22}x_2^2 = 0, \quad (3.10)$$

where $a_{ij}$ are integers.

It is known that the equation (3.10) is reduced to the diagonal form by rotatiing through $\varphi$, where:

$$tg\varphi = (a_{22} - a_{11} \pm \sqrt{(a_{22} - a_{11})^2 + 4a_{12}^2})/2a_{12}. \quad (3.11)$$

If $a_{11} = a_{22}$, then taking account of (3.11) we get $tg\varphi = \pm 1$ and $\varphi = \pi/4$ or $\varphi = -\pi/4$. In these cases the rotation matrix C respectively is:

$$C = \begin{pmatrix} \sqrt{2}/2 & -\sqrt{2}/2 \\ \sqrt{2}/2 & \sqrt{2}/2 \end{pmatrix} \text{ or } C = \begin{pmatrix} \sqrt{2}/2 & \sqrt{2}/2 \\ -\sqrt{2}/2 & \sqrt{2}/2 \end{pmatrix}, \quad (3.12)$$

and it contains the same irrationalities.

Based on (3.5) with the homothetic with $k = \sqrt{2}$ we get the integer deformation matrices from the rotation matrix (3.12):

$$C_d = \begin{pmatrix} 1 & -1 \\ 1 & 1 \end{pmatrix} \text{ or } C_d = \begin{pmatrix} 1 & 1 \\ -1 & 1 \end{pmatrix}. \quad (3.13)$$

Since $a_{11} = a_{22}$ we get diagonal Diophantine equations after the deformation conversion (3.13):

$$(a_{11}+a_{12})x'^2_1+(a_{11}-a_{12})x'^2_2=0, \ (a_{11}-a_{12})x'^2_1+(a_{11}+a_{12})x'^2_2=0.$$

If $\varphi = \pi/3$ the rotation matrix has the form:

$$C = \begin{pmatrix} 1/2 & -\sqrt{3}/2 \\ \sqrt{3}/2 & 1/2 \end{pmatrix}. \qquad (3.14)$$

Matrix (3.14) corresponds to the fourth case in which $k$ does not exist.

Let's consider the general algebraic Diophantine equation of order $m$ (3.1).

If the equation (3.1) corresponds to the central hypersurface, then assume that it can be reduced to the canonical form by transformation of the motion:

$$F(x'_1,...,x'_n) = \sum_{i=1}^{n} a'_{ii}(x'_i)^m + a'_0 = 0, \qquad (3.15)$$

where the coefficients, in general, are not integers. It is performed for algebraic Diophantine equation of the second order without assuming.

If the equation (3.1) does not correspond to the central hypersurface, then assume that it can be reduced to the canonical form by motion transformation:

$$F(x'_1,...,x'_n) = \sum_{i=1}^{n-1} a'_{ii}(x'_i)^m + a'_n x'_n = 0, \qquad (3.16)$$

where the coefficients in general are not integers. It is performed for algebraic Diophantine equation of the second order without assuming.

Conversion movement which reduces equation (3.1) to the canonical form (3.15) or (3.16) is described by (3.2), where the coefficients $c_{ij}$ and $c_i$ in general are not integer. If $c_i$ is integer, then you can do it all in the above mentioned first and second cases with integer coefficients $c_{ij}k$ by a homothetic transformation $x'_i = kx''_i$ and thus make the transformation (3.3) completely integer. Then the equation (3.15) can be transformed to diagonal Thue equation, which was discussed above:

$$F(x''_1,...,x''_n) = \sum_{i=1}^{n} k^m a'_{ii}(x''_i)^m + a'_0 = 0, \qquad (3.17)$$

and the equation (3.16) can be transformed to the following form:

$$F(x"_1,...,x"_n) = \sum_{i=1}^{n-1} k^{m-1} a'_{ii}(x"_i)^m + a'_n x"_n = 0. \qquad (3.18)$$

We use the fact that any algebraic Diophantine equation of order k corresponds to an algebraic hypersurface of the same order with integer coefficients. If perform integer affine transformation for the algebraic hypersurface, it is also transformed into an algebraic hypersurface of order k with integer coefficients, and this hypersurface will correspond to an algebraic Diophantine equation of the same order. Therefore, all the coefficients of (3.17) and (3.18) are integer.

Note. Diophantine homogeneous equation and Thue equations correspond to the central hypersurfaces which center coincide with the coordinate origin. Therefore, the transformation of these equations to a diagonal form does not include the transfer and the first condition, i.e. the transfer of integer coordinates, is done automatically.

## 4. ESTIMATION OF THE NUMBER OF INTEGER (NATURAL) SOLUTIOMS OF SOME TYPES OF NONDIAGONAL ALGEBRAIC DIOPHANTINE EQUATIONS

The above considered deformation transformation in general does not lead algebraic Diophantine equations of the order higher than second order to the diagonal form. Therefore, it is necessary to examine estimations of the number of integer (natural) solutions for **nondiagonal** algebraic Diophantine equations with values $m \geq 3$.

It is known [14] that the number of integer solutions for not only diagonal algebraic Diophantine equation of the order $m \geq 3$ and with variables $n > 2$ has the following upper bound in the hypercube with the sides $[-N, N]$:

$$R^m_n(N) \leq m(2N+1)^{n-1}, \qquad (4.1)$$

where Diophantine algebraic equation in general can be reduced.

Pila formula (2.12) can be used also for an irreducible not only diagonal algebraic equation of the order $m \geq 3$ with variables $n > 2$ in the hypercube with the sides $[-N, N]$.

Assertion 1. Let's suppose an algebraic Diophantine equation $F(x_1,...,x_n) = 0$ of not only diagonal form satisfying the condition $F(0,...0, x_l, 0,...,0, x_k, 0,...,0) = 0$ that is an irreducible algebraic Diophantine equation of the kind $g \geq 1$. Then we have the following upper bound for the number of integer solutions of algebraic Diophantine equation in the hypercube with the sides $[-N, N]$:

$$R_n(N) << N^{n-2}. \tag{4.2}$$

Proof. Diophantine equation $F(0,...0, x_l, 0,...,0, x_k, 0,...,0) = 0$ has a finite number of solutions according to the Siegel theorem [1]. The maximum number of integer values that are taken by other variables is equal to $(2N+1)^{n-2}$, so the relation that corresponds to (4.2) is true:

$$R_n(N) \leq R_2(N)(2N+1)^{n-2} = O(1)O(N^{n-2}) = O(N^{n-2}).$$

Suppose a Diophantine Thue equation of the order $m$:

$$F_m(x_1,...,x_{n-1}) + c_n x_n^m + F_0 = 0, \tag{4.3}$$

where $F_m(x_1,...,x_{n-1}) = 0$ is the nondiagonal homogeneous algebraic Diophantine equation, and $c_n, F_0$ are integers.

Assertion 2. If $x_n = \sqrt[m]{-F_0/c_n}$ is an integer, then the equation (4.3) has at least one integer solution.

The proof follows from the fact that the homogeneous algebraic Diophantine equation always has the trivial solution.

Assertion 3. If $x_n = \sqrt[m]{-F_0/c_n}$ is an integer and the equation $F_m(x_1,...,x_{n-1}) = 0$ has at least one nontrivial solution $(x_{10},...,x_{l0})$, then the equation (4.3) has an infinite number of integer solutions and there is at least as much as $(2N+1)/\max\{x_{10},...,x_{l0}\}$ integer solutions in the hypercube with the sides $[-N, N]$.

Proof. The equation $F_m(x_1,...,x_{n-1}) = 0$ is a homogeneous equation, so if it has at least one nontrivial integer solution, then it has an infinite number of integer solutions and there is at least as much as $(2N+1)/\max\{x_{10},...,x_{l0}\}$ integer solutions in the hypercube with the sides $[-N, N]$. If $x_n = \sqrt[m]{-F_0/c_n}$ is an integer, then the equation (4.3) has an infinite number of integer solutions and there is at least as much as $(2N+1)/\max\{x_{10},...,x_{l0}\}$ integer solutions in the hypercube with the sides $[-N, N]$.

## 5. CONCLUSION AND SUGGESTIONS FOR FURTHER WORK

The next article will continue to study integer conversions that maintain the asymptotic behavior of the number of integer solutions for the algebraic Diophantine equation.

## 6. ACKNOWLEDGEMENTS

Thanks to everyone who has contributed to the discussion of this paper.